\input amstex
\documentstyle{amsppt}
\magnification=\magstep0
\define\cc{\Bbb C}

\define\r{\Bbb R}

\define\N{\Bbb N}

\define\jj{\Bbb J}

\define\A{\Cal A}
\define\h{\Cal D}
\define\E{\Cal E}

\define\f{\Cal S}

\define\la{\lambda}
\define\om{\omega}

\define\e{\varepsilon}
\define\va{\varphi }
\define\CB#1{\Cal C_b(#1)}
\define\st{\subset }

\topmatter
 \title
 Comparison of spectra of  absolutely regular   distributions and applications.
  \endtitle
 \author
  Bolis Basit and  Alan J. Pryde
\endauthor
 \abstract
{We study  the reduced Beurling spectra $sp_{\Cal {A},V} (F)$  of functions
  $F \in L^1_{loc}
 (\jj,X)$ relative to certain function spaces $\Cal{A}\st L^{\infty}(\jj,X)$ and $V\st L^1 (\r)$ and compare them with other spectra including the weak Laplace spectrum. Here $\jj$ is $\r_+$ or $\r$ and $X$ is a Banach space. If  $F$  belongs to  the space $ \f'_{ar}(\jj,X)$ of absolutely regular distributions and has uniformly continuous indefinite integral with $0\not\in sp_{\A,\f(\r)} (F)$ (for example if F is slowly oscillating and $\A$ is $\{0\}$ or $C_0 (\jj,X)$), then $F$ is ergodic.   If $F\in  \f'_{ar}(\r,X)$ and $M_h F (\cdot)= \int_0^h F(\cdot+s)\,ds$ is bounded for all $h > 0$ (for example if $F$ is ergodic) and if $sp_{C_0(\r,X),\f} (F)=\emptyset$, then ${F}*\psi \in C_0(\r,X)$ for all $\psi\in \f(\r)$.
  We show that tauberian theorems for Laplace
transforms follow
   from results about  reduced spectra.    Our results are more general than previous ones and we demonstrate this  through  examples}
\endabstract
\endtopmatter
\rightheadtext{  reduced spectra}
\leftheadtext{B. Basit,  A. J. Pryde}
  \TagsOnRight
\document
\pageno=1 \baselineskip=18pt

26 August 2011

\head{\S 1. Introduction}\endhead

The goal of this paper \footnote {AMS subject classification 2010:  Primary  {47A10, 44A10} Secondary {47A35, 43A60}.
\newline\indent Key words: Reduced Beurling,  Carleman, Laplace and weak Laplace  spectra, almost periodic, asymptotically almost periodic}
  is to study the asymptotic behaviour of certain locally integrable functions $F: \jj\to X $ where $\jj$ denotes $\r$ or $\r_+$ and $X$ is  a complex Banach space. Such a study has a long history.  It is motivated by  Loomis's theorem (see [25] and [24, Theorem 4, p. 92, p. 97]) which gives spectral conditions under which a function $F: \r\to \cc$ is almost periodic and by the tauberian theorem of Ingham  (see [22] and [2, Theorem 4.9.5, p. 326]) which gives conditions under which $\lim_{t\to \infty} F(t)=0$. Many notions of the spectrum of a function have since been introduced in order to obtain (vector valued) extensions of these results and we will review and compare  some of these in this paper. In particular
   we develop
  the reduced  spectrum $sp_{\A} (F)$ of $F$ relative
 to   various closed subspaces $\A$ of $BUC(\jj,X)$, a spectrum   that was introduced  before in this context  (see [24, Chapter 6.4, p. 91], [6], [7], [2, p. 371], [8] and [17]). Typically the reduced spectrum gives stronger results than other spectra as we shall see.
  The style of  definition and its properties are similar to  those of  the Beurling spectrum (see [6, Theorem 4.1.4]) and it
 is widely applicable.  For $ F\in BUC(\jj,X)$ and $\A \st BUC(\jj,X)$, there is  also an operator theoretical approach using $C_0$-semigroups and groups (see [19] for example).
  Attempts by  Minh [26] to define a reduced spectrum of bounded not necessarily uniformly continuous functions $F\in BC(J,X)$ using the operator theoretical approach  failed  even for the case $\jj=\r$ (see [27]).

   The reduced spectrum was replaced by the local Arveson spectrum   in [2] and  [5, p. 296] and by the Laplace spectrum  $sp^{\Cal{L}} (F)$ in [3], [14] and [15]. To extend the theory, a smaller spectrum, the half-line spectrum $sp_+(F)$, was introduced in [4, p. 474]. Then Chill introduced a still smaller spectrum, the weak Laplace spectrum $sp^{w\Cal{L}} (F)$ in [17, p. 25] and [18, Definition 1.1]. Typically $sp^{w\Cal{L}} (F) \st sp_+(F) \st sp^{\Cal{L}} (F)$.

    The important  Theorem 5 of Chill and Fasangova [19] (see also [10, Theorem 3.10]), establishing  that
   the reduced spectrum coincides with the local Arveson spectrum for the group  induced by the shifts on a quotient space,   shows that some results of
  Arendt and Batty in [2] and [3] using the local Arveson  spectrum  follow from  earlier results
  in [6] and [7]. This point is noted in  [19, Theorem 10].  Similarly, noting that  $sp_{\A} (F) \st  sp_{C_0} (F) \st sp^{w\Cal{L}} (F)\st  sp^{\Cal{L}} (F)$ (Proposition 4.2), the general tauberian  theorems [3, Theorem 2.3] and  [15, Theorem 4.1] using
   the the Laplace spectrum $sp^{\Cal{L}} (F)$ are  consequences of Theorems  4.2.5 and 4.2.6  of [6] which use the reduced spectrum. The latter  theorems are stronger than the main tauberian theorem of [15] and their proofs are simpler. A difficulty with $sp^{\Cal{L}} (F)$ and $sp^{w\Cal{L}} (F)$ is that it is unclear whether they satisfy the useful property  $sp^{*} (F*g)\st sp^{*} (F)\cap $ supp $\,\,\widehat {g}$, for appropriate  $g$, even when $F$ is bounded; but see Proposition 3.4 (i) for reduced spectra.

    In (3.4)  we consider a more general spectrum $sp_{\A,V} (F)$, the reduced spectrum of $F\in L^1_{loc} (\jj,X)$ relative to $(\A,V)$, where $V\st L^1(\r)$, a spectrum  closely related to the  one  defined in (3.4$^*$) which was first studied in [10].    We are able to strengthen further the improvements  made by Chill ([17, Lemma 1.16] and [18, Proposition 1.3, Theorem 1.5,  Corollary 1.7]). In particular  we replace   $sp^{w\Cal{L}} (F)$  by  $sp_{\A,V} (F)$ for some $\A \supset C_0 (\r_+,X)$. We are able to consider functions whose Fourier transforms  are not regular distributions  (see Example 3.12) and  avoid some geometrical restrictions  on  $X$  that were imposed in [17, Theorem 1.23] and [18, Proposition 2.1] for example.
     Moreover, our methods   lead  to  new results  for the weak Laplace spectrum (see Theorem 4.3). Finally, spectral criteria  for solutions of evolution equations  on $\r_+$ or $\r$ are readily  related  to  reduced spectra  (see  Theorem 4.6 and [7, (1.7)]).

In section 2 we describe our
notation  and prove some preliminary results. We are particularly interested in functions $F\in \f'_{ar} (\jj,X)$, the space of absolutely regular distributions (see (2.1)).

 In section 3 we study  $sp_{\A,V} (F)$, where  $V$ is one of the spaces $\h=\h(\r)$, $\f=\f(\r)$ or $L^1=L^1(\r)$. If  $F\in L^{\infty}(\r,X)$ and  $\A= 0:= \{0\}$, then by [10, (3.3)] (see also [5, Proposition 4.8.4, (4.26)], [28, 0.5,  p. 19], [29,  p. 183])

(1.1)\qquad  $sp_{0,\f} (F)= sp_{0,L^1} (F)=sp^{\Cal {C}} (F)=sp^{B} (F)=sp_0 (F)$,

\noindent where $sp^{\Cal {C}} (F)$ is the Carleman spectrum and  $sp^{B} (F)$ is the  Beurling spectrum.
In  Proposition 3.1 (i) we prove that our definition (3.4) coincides  with (3.4$^*$) for the spaces $V\in \{ \h,\f,L^1\}$.  In  (3.3) and  Proposition 3.2  we study the conditions imposed on $\A$ and relate them to others in the literature, in particular to the translation-biinvariance  used in  [17, Definition 1.2, p. 17] and [3, \S2].  In Remark 3.3 we show that the converse of Proposition 3.2 (i) is false and that $F\in \A$  does not imply in general that $sp_{\A} (F)=\emptyset$. We develop  some basic properties of the reduced spectrum in Proposition 3.4. Our main results are stated in Theorems  3.5, 3.6, 3.7 and 3.8.

In  Theorems  3.5 and 3.6 we prove  ergodicity results for functions  $F\in \f'_{ar} (\jj,X)$. If $\jj=\r$, $0\not\in sp_{0,\f} (F)$ and  the indefinite integral $PF\in UC(\r,X)$, then  $PF\in BUC(\r,X)$ and $F$ is ergodic. For a variety of classes $\A$ including $C_0 (\jj,X)$, $EAP_0(\jj,X)$,  $EAP(\jj,X)$, $AAP(\jj,X)$ and  $AP(\r,X)$, if $0\not\in sp_{\A,\f} (F)$ then $F$ is ergodic.

Theorem 3.7 deals with functions $F \in\f'_{ar}(\jj,X)$ with $sp_{\A,\f} (F)$  countable. It is a generalized tauberian theorem providing spectral conditions under which $F$ has various types of  asymptotic behaviour.  For example (Theorem 3.7 (v)), if $M_h F$ is bounded for each $h >0$, $sp_{\A,\f} (F)$ is countable and $\gamma_{-\om} F$ is ergodic for each  $\om\in sp_{\A,\f} (F)$, then $(\overline {F}*\psi)|\, \jj \in \A$ for all $\psi\in \f (\r)$. (For the definition of $\overline {F}$ see (2.5)).  Versions of Theorem 3.6,  Theorem 3.7 (i)- iv) and Corollary 3.9   are already known when $\jj=\r_+$ and $sp_{\A,\f} (F)$ is replaced by the larger spectrum $sp^{w\Cal{L}} (F)$ (see  Remark 4.4 (i)).
      Corollary 3.10  states that if $F\in L^1_{loc} (\jj,X)$  with $sp_{C_0(\jj,X),\h} (F)=\emptyset$ and if the convolution $(\overline {F}*\psi)|\,\jj$ is  uniformly continuous for some $\psi\in \h(\r)$
 then $\overline {F}*\psi \in C_0(\r,X)$. For the case $\jj=\r$,  Chill [18, Proposition 3.1]  obtained this same conclusion under the stronger assumptions that $F\in L^1_{loc}(\r,X)$ and $\widehat{F}\in \f'_{ar}(\r,X)$. In particular, if $F\in L^p (\r,X)$ where $1\le p < \infty$, then $F$ satisfies the assumptions of Corollary 3.10. However, as is well-known, when $p > 2$ there are functions $F\in L^p (\r,X)$ for which $\widehat{F}$ is not a regular distribution and so the result of [18] does not apply. Even when $1\le p\le 2$ special geometry on $X$ is required in order that every $F\in L^p (\r,X)$ has a Fourier transform which is regular.

In Proposition 4.1  we  establish some  properties  shared by the  weak Laplace, Laplace
and Carleman
spectra. In  Proposition 4.2  we prove the inclusion \,\,$sp_{\Cal {A},V} (F)$
 $\subset sp^{w\Cal{L}}(F)$  for $F\in \f'_{ar} (\r_+,X)$ and ${\Cal{ A}\supset
C_0 ({\Bbb R}_+,X)}$. This enables us to prove a new tauberian theorem (Theorem 4.3) and also to deduce the main tauberian results of Chill mentioned  above (see  Remark 4.4). In Example 4.5 we demonstrate how to use Proposition 4.1 to calculate Laplace spectra. Finally,
 in Theorem 4.6 we obtain a spectral condition satisfied by  bounded mild solutions
of the evolution equation $\frac{d u(t)}{dt}= A u(t) +  \phi (t) $, $u(0)\in{X}$, $t\in  {\jj}$,
  where $A$ is a closed linear operator   on ${X}$ and $\phi\in L^{\infty} (\jj, {X})$.  This  generalizes earlier results  where it is assumed that $u,\phi \in BUC(\jj,X)$ (see [5, Proposition 5.6.7] and [7, Theorem 3.3, Corollary 3.4]).

 \head{\S 2. Notation, Definitions and preliminaries}\endhead

Throughout the paper $\r_+=[0,\infty)$, $\r_-=(-\infty,0]$, $\jj\in\{\r_+,\r\}$, $\N= \{1,
2,\cdots \}$, $\cc_+
 =\{\la\in \cc: \text{Re \,}\la
>0\}$ and $\cc_-=\{\la\in \cc: \text{Re \,}\la < 0\}$.
  By $X$ we  denote a   complex Banach space. If $Y$ and $Z$ are locally
 convex topological spaces, $L(Y,Z)$  denotes the space of all bounded linear operators from $Y$ to
 $Z$.
   The Schwartz spaces of test functions and rapidly decreasing functions are denoted by $\h(\r)$ and $\f
   (\r)$ respectively.  Then  $\h' (\r,X)= L(\h(\r),X)$ is the  space of $X$-valued distributions and   $\f' (\r,X)= L(\f(\r),X)$ is the  space of $X$-valued tempered distributions
   (see [5, p. 482],[32, p. 149] for $X=\cc$). The space  of absolutely regular distributions is defined by

(2.1) \qquad $\f'_{ar} (\jj,X)=\{ H\in L^1 _{loc}(\jj,X): H\va \in L^1(\jj,X)$ for all $\va\in \f(\r) \}$.

 \noindent   The action of an element $ S\in \h' (\r,X)$ or  $\f' (\r,X)$ on
  $\va\in\h(\r) $ or $\f(\r)$ is denoted  by $<S,\va>$.
 If $F$ is an  $ X$-valued function defined on  $\jj$ and $s\in \jj$
then   $F_s$,  $\Delta _sF$,  $|F|$  stand for the functions defined on
$\jj$ by $F_{s}(t) = F(t+s)$, $\Delta _sF (t)= F_{s}(t)-F (t) $ and $|F| (t)= ||F(t)||$.
Also  $||F|| _{\infty} =
  \text {  sup}_{t\in \jj} ||F(t)||$.
 If $F \in L_{loc}^1 (\jj, X)$ and $h>0$, then $PF$,  $M_hF$  and $\check {F}$ (when $\jj=\r$)  denote  the
$indefinite\,\, integral$,  $mollifier$  and $reflection$ of $F$ defined
respectively by $PF(t) = \int_{0}^{t} F (s)\,ds$, $M_h F (t)=
(1/h)\int_0^h F (t+s)\,ds$ for $t\in \jj$ and $\check {F} (t) =F(-t)$ for $t\in \r$.
 For $g\in
   L^1(\r)$ and $F\in
   L^{\infty}(\r,X)$ or $g\in
   L^1(\r,X)$ and $F\in
   L^{\infty}(\r)$ the $Fourier\,\, transform$ $\widehat{g}$ and $convolution$ $F*g$ are
   defined respectively by $\widehat{g} (\om)=\int_{-\infty}^{\infty} \gamma_{-\om} (t)\, g(t)\,
   dt$ and $F*g (t)= \int_{-\infty}^{\infty}F(t-s) g(s)\,
   ds$, where  $\gamma_{\om} (t)= e^{i\, \om t} $ for $\om\in\r$.
    The $Fourier\, transform$ of $ H\in \f' (\r,X)$ is
  the tempered distribution $\widehat{H}$  defined by

  (2.2)\qquad  $<\widehat{H},\va> =<H,\widehat {\va}>$\,\, for all $\va\in
   \f(\r)$.

\noindent Set $\widehat{\h}(\r)=\{\widehat{\va}: \,\,\, \va\in \h
(\r)\} \st \f(\r)$.
    The $Fourier\, transform$ of $F \in L_{loc}^1 (\r, X)$
   is the distribution $\widehat {F}\in L(\widehat{\h}(\r),X)$
   defined by

  (2.3)\qquad $<\widehat{F},\psi> =<F,\widehat{\psi}>$\,\, \, for all $\psi\in
  \widehat{ \h}(\r)$.

\noindent Throughout the paper all integrals are Lebesgue-Bochner
integrals (see [5, pp. 6], [20, p. 318],
   [21, p. 76]). All convolutions  are understood as convolutions of functions defined on $\r$. Given $F\in  W (\jj,X)$ where

 (2.4) \qquad     $ W (\jj,X) \in  \{L^1_{loc} (\jj,X), \f'_{ar} (\jj,X), L^{\infty} (\jj,X)\}$,

\noindent     we denote by $\overline{F}: \r\to X$ the function given  by

 (2.5) \qquad    $\overline{F}|\jj =F$ and,  if $\jj=\r_+$, $\overline{F}|(-\infty,0) =0$.

\noindent Then $\overline{F}\in W(\r,X)$. In addition, if $g\in L^{\infty}_c (\r)=\{f\in L^{\infty} (\r): f $ has compact support $\}$, then for some constant $t_{g}$

 (2.6) \qquad $\overline{F}*g\in W (\r,X) \cap C(\r,X)$ and,  if $\jj=\r_+$,  $\overline{F}*g (t)=0$ for all $t\le t_{g}$.

\noindent  It follows that if  $h > 0$ and  $s_h= (1/h)\chi_{(-h,0)}$, where   $\chi_{(-h,0)}$ is the characteristic function of $(-h,0)$, then

(2.7) \qquad  $\overline{F}*s_h\in W(\r,X)\cap C(\r,X)$, \,\,
 $M_h F= (\overline{F}*s_h) |\jj$ and,

\qquad\qquad\,\,\, if $\jj=\r_+$,  $ \overline{F} *s_h (t) =0$ for all $t \le -h$.

 We  use convolutions of functions $F\in W=W(\jj,X)$ and $g \in V= V(\r)\in\{ \h(\r), \f(\r), L^1(\r)\}$, with

(2.8) \qquad $V=\h(\r)$ if $W= L^1_{loc} (\jj,X)$, $V=\f(\r)$ if $W=\f'_{ar}(\jj,X)$ and

\qquad \qquad \,\,\,  $V=L^1(\r)$ if $W=L^{\infty}(\jj,X)$.

 \noindent The following properties of the convolution  are repeatedly used
(see [32, p. 156, (4)] and [30, 7.19 Theorem (a), (b), pp.179-180] when $X=\cc$):

\noindent If $ F\in W(\jj,X)$ and $\va\in V(\r) $   where  (2.4) and (2.8)  are satisfied, then

(2.9)\qquad  $\overline{F}*\va\in W(\r,X) \cap C(\r,X)$.

\noindent Indeed, the cases $W=L^1_{loc}(\jj,X)$ and   $W=L^{\infty}(\jj,X)$ are obvious. If   $F\in \f'_{ar}(\jj,X)$, then
$|F|\in \f'_{ar} (\jj,\cc)$. By [23, Theorem (b)] there is an integer  $k \in \N$  such that

(2.10.) \qquad $\frac{|F|}{w_k}= f\in L^1(\jj)$, where $w_k (t)= (1+t^2)^k$.

\noindent  Using (2.10), we easily conclude (2.9).

 Moreover, if $\psi \in V(\r)$ or $\psi\in L_c^{\infty}(\r)$, then

(2.11)\qquad  $(\overline{F}*\va)*\psi= \overline{F}*(\va*\psi)=(\overline{F}*\psi)*\va$.

Now let $F\in \f'_{ar} (\r,X)$ and $\va\in\f(\r)$ be such that $\widehat{\va}\in \h(\r)$ and  $\widehat{\va}=1$ on $[-\delta,\delta]$ for some $\delta > 0$. Then

(2.12)\qquad $0\not \in sp_{0,\f} (F-F*\va)$.

\noindent Indeed, if $\psi\in\f(\r)$, supp $\widehat{\psi}\st [-\delta,\delta]$ and  $\widehat{\psi} (0)=1$, then $\va*\psi =\psi$. So by (2.11), $(F-F*\va)*\psi=0$.

For the benefit of the reader we include the proofs of the following elementary  but necessary results.

\proclaim{Lemma 2.1}  Let $f\in L^1 (\r)$  with $\widehat {f} \not =
0$ on a compact set $ K$. Then  there exists $g\in L^1 (\r)$ such that
$\widehat {g} \cdot \widehat {f} =1$ on $ K$. Moreover, one can choose $g$
such that $\widehat{g}$ has compact support and, if $f\in \f (\r)$,
with $g\in \f (\r)$.
\endproclaim

\demo{Proof} Choose a bounded open set $U$ such that $K\st U$ and
$\widehat {f} \not = 0$ on $\overline {U}$ the closure of $ {U}$. By
[16, Proposition 1.1.5 (b), p. 22], there is $k \in L^1(\r)$ such
that $\widehat {k} \cdot \widehat {f} =1$ on $ \overline {U}$. Now, choose
$\va\in \h(\r)$ such that $\va=1$ on $K$ and supp $\va\st
\overline {U}$.  Also choose $ \psi\in
\f(\r)$ such that $\widehat{\psi}=\va$ and take  $g= k*\psi$. Then
$\widehat {g}$ has compact support and if $f\in \f(\r)$, then
$\widehat{g}\in \h(\r)$ and so $g\in \f(\r)$. $\square$
\enddemo

In the following  lemma  $\psi$ will denote an element  of $ \f
(\r)$ with the  properties:

\qquad $\widehat {\psi}$  has compact support,  $\widehat
{\psi} (0)=1$ and $\psi $ is non-negative.

\noindent An example of such  $\psi$ is given by
$\psi=\widehat{\va}^2$, where $\va(t)=a\, e^{\frac{1}{t^2-1}}$ for
$|t|\le 1$ and  $\va=0$ elsewhere on $\r$ for  some suitable
constant $a$.

\proclaim{Lemma 2.2} (i)
 The sequence $\psi_n(t)= n\, \psi (n\,t)$ is an approximate identity
for the space of uniformly continuous functions $UC(\r,X)$, that
is $\lim _{n\to \infty}||u*\psi_n- u||_{\infty}= 0$  for all $u\in UC(\r,X)$.

(ii) $\lim_{h\searrow 0}||M_h u- u||_{\infty} = 0$  for all $u\in UC(\jj,X)$. In particular if $M_h u\in BUC(\jj,X)$ for all $h > 0$ then $ u\in BUC(\jj,X)$.
\endproclaim

\demo {Proof} (i) Given  $u\in UC(\r,X)$ and $\e >0$ there exists
$k >0$ such that  $||u (t+s)-u(t)|| \le k  |s|+\e$ for all
$t,s\in \r$. In particular $u\in \f'_{ar} (\r,X)$.  Also,
 $u*\psi_n (t)- u(t)= \int_{-\infty}^{\infty} [u (t-\frac{s}{n})-u(t)]
 \psi (s)\,ds$ which gives $||u*\psi_n - u||_{\infty}\le (k/n)\int_{-\infty}^{\infty}
 |s|\psi (s)\,ds+ \e \int_{-\infty}^{\infty}
 \psi (s)\,ds$ and  (i) follows.

 (ii) Since $||M_h u- u||_{\infty}\le \text{\, sup} _{t\in \jj, 0 \le s \le h}||u (t+s)- u(t)||$, part (ii) follows.
  \P
\enddemo
\proclaim{Lemma 2.3}  Let $F\in W (\r,X)$  and $g\in V(\r)$ with   (2.4) and (2.8) satisfied.

(i) If $F|\r_- =0$, then
 $(F*g) |\r_- \in C_0 (\r_-,X)$.

(ii) If $F|\r_+ =0$, then
 $(F*g) |\r_+ \in C_0 (\r_+,X)$.
\endproclaim

\demo {Proof} (i)  The cases $W= L^1_{loc} (\r,X)$ and $W= L^{\infty} (\r,X)$  can be shown by simple calculations. If $F\in \f'_{ar} (\r,X)$, then from (2.10) ${|F|}/{w_k}= f\in L^1(\r)$ for some $w_k (t)= (1+t^2)^k$. Since $\va\in \f(\r)$,
 $||w_k \va||_{\infty}=c_k < \infty $.  It follows that
 $||F*\va (t)||= ||\int_0^{\infty} \va (t-s) F(s)\, ds|| \le  \int_0^{\infty} |\va| (t-s) |F|(s)\, ds \le  c_k \int_0^{\infty} \frac{w_k (s)}{w_k (t-s)} f(s)\, ds$. Since $\frac{w_k (s)}{w_k (t-s)}\le 1$ for each $t\le 0, s \ge 0$ and $\lim_{t \to -\infty} \frac{w_k (s)}{w_k (t-s)}$ $=0$
for each $s \ge 0$, it follows that $\lim_{t\to -\infty} ||{F}*\va (t)||=0$ by the Lebesgue convergence theorem. By (2.9) the result follows.

 (ii) This follows by applying part (i) to $\check {F}$.
 \P
\enddemo
\head{\S 3. Reduced spectra for regular distributions}\endhead

In this section we introduce  the reduced spectrum $sp_{\A,V} (F)$ of a function $F\in L^1_{loc} (\jj,X)$  relative to $\A, V$, where $\A \st L^{\infty} (\jj,X)$ and $V\st L^1(\r)$. We usually impose the following conditions on $\A$.

(3.1) \qquad  $\A$ is a  closed subspace of $
L^{\infty}(\jj,X)$ and is $BUC$-$invariant$; that is

\qquad\qquad \,\,\, if $F\in BUC(\r,X)$ and $F|\,\jj\in \A$,
then $F_t|\,\jj\in \A$ for each $t\in \r$.

The property of being $BUC$-$invariant$ was first introduced in  [6, $(P.\Lambda)$, Definition 1.3.1] and called the Loomis property  for classes $\A \st BUC(\jj,X)$. The notion was extended to classes $\A\st L^1_{loc} (\jj,X)$ in  [8, (1.III$_{ub}$)]. In [10], this property was called  $C_{ub}$-invariance. In the proof of Theorem 2.2.4 in [6, p. 13], it is shown that if $\jj=\r_+$  and $\A$ satisfies  (3.1), then

(3.2) \qquad  $ C_0(\r_+,X) \st \A_{ub}$, where  $\A_{ub}=\A \cap BUC(\r_+,X)$.

We note that if $\jj=\r$, then  $\A$ is $BUC$-invariant if and only if $\A \cap BUC (\r,X)$ is a translation invariant subspace of $BUC (\r,X)$.
If $\jj=\r_+$, then $\A$ is  $BUC$-invariant if and only if $\A \cap BUC(\r_+,X)$ is  a positive invariant subspace of $BUC (\r_+,X)$ (that is $F_t \in \A$ whenever $F\in \A$, $t\ge 0$) with the additional  property that $F\in \A$ whenever $F\in BUC (\r_+,X)$ and  $ F_t \in \A $ for some  $t \ge 0$.
  Such subspaces of $BUC(\r_+,X)$ were  called
 $translation$-$biinvariant$
   ([17, (1.1), p. 17], [3]). So, for a closed subspace  $\A $ of $BUC(\r_+,X)$,

  (3.3) \qquad   $\A $  is $BUC$-invariant if and only if $\A$ is translation-biinvariant.

For $\A$ satisfying  (3.1),  $V\st L^1 (\r)$ and $F\in L^1_{loc} (\jj,X)$,  a point $\om\in \r$ is called $(\A, V)$-$regular$
for $ F$ or $\overline{F}$, if there is $\va\in V$ such that $\widehat {\va}(\om)\not =0$ and
$(\overline{F}*\va)|\,\jj \in \A$. The $reduced$ $Beurling$ $spectrum$ of $F$ or $\overline{F}$
relative to $(\A,V)$ is defined by

(3.4) \qquad $sp_{\A,V} (F) =\{\om\in \r: \om$ is not an
$(\A,V)$-regular point for $F \}=$

\qquad\qquad  \,\, $ \{\om\in \r: \va \in V, (\overline{F}*\va)|\, \jj\in \A
$ implies  $\widehat{\va} (\om)=0 \}=sp_{\A,V} (\overline{F})$,

\noindent provided  the convolution $\overline{F}*\va$ and the restriction $\overline{F}*\va)|\, \jj$ are defined for all $\va\in V$. Clearly, $sp_{\A,V} (F)$ is a closed subset of $\r$.
Further, if $F\in L ^1_{loc} (\r,X)$ and $\A \st L^{\infty} (\r_+,X)$  we also explore  the following spectrum first introduced  in [10, Definition 3.1]

(3.4$^{*}$) \qquad $sp_{\A,V} (F)=\{\om\in \r: \va \in V, (F*\va)|\, \r_+\in \A
$ implies  $\widehat{\va} (\om)=0 \}$.

\noindent We give conditions in Proposition 3.1 under which $sp_{\A,V} (F)= sp_{\A,V} (F|\r_+)$.

If $F\in \f'_{ar} (\r_+,X)$, then $sp_{\A,\f} (F)\st sp_{\A,\h} (F)$  since $\h (\r) \st \f(\r)$. Moreover, the  first inclusion might be proper. For example take  $F\in L^{\infty}(\r,X)$ with  $sp^B(F)=sp_{0,\f} (F) $  (see (1.1))   uncountable but not  $\r$.    If $0\not =\va\in \h(\r)$, then  $F*\va \not =0 $ otherwise $sp^B(F)$ is  countable. It follows that  $sp_{0,\h} (F) =\r$.

For  $ F\in L^{\infty}(\jj,X)$ and $V= L^1(\r)= L^1$ we write $sp_{\A} (F)= sp_{\A,L^1} (F)$.

 If $ F\in W(\jj,X)$  and  $V=V(\r)$  satisfies  (2.8), then the convolution $\overline{F}*g$ and the restriction $(\overline{F}*g)|\, \jj$ are defined for all $g\in V(\r)$. So, $sp_{\A,V} (F)$    is well defined.
  This is an
extension of the definitions in  [6, (4.1.1)], [7, (2.9)] and [17, Definition 1.14, p. 24]. In those references the conditions on $\A$
 are more restrictive and $ F\in L^{\infty}(\r,X)$.

 If  $ F\in
\f'_{ar}(\r,X)$, then
$sp_{0,\f} (F)= $ supp
$\widehat{F} = sp^{\Cal{C}} (F) $,
where $sp^{\Cal{C}} (F) $ is the  the Carleman spectrum. See [28, Proposition 0.5, p. 22].

   If $ F\in L^{\infty}(\jj,X)$, then $\overline{F}*f\in BUC(\r,X)$ for all $f\in L^1(\r)$. It follows that \,\,\,

(3.5)\qquad  $sp_{\A} (F)= sp_{\A_{ub}} (F)$, where  $\A_{ub}=\A\cap BUC(\jj,X)$.

\noindent  A sufficient condition to have the property  $sp_{\A} (F)=\emptyset$   for each $F\in \A \st L^{\infty} (\jj,X)$ is the following

  (3.6)\qquad  $(\overline{F}*f)|\, \jj\in \A_{ub}$ for each $F\in \A$ and $ f\in L^1(\r)$.

 Examples of spaces  $\A$ satisfying  (3.6) include (using $\A(\jj,X)=\A(\r,X)|\jj $)

\smallskip

\noindent  $\{0\}$,\quad  $C_0 (\jj,X)$,\quad $AP= AP(\r,X)$,\quad $LAP_b (\r,X)$,\quad
  $AA= AA(\r,X)$,\quad

\noindent $ EAP(\jj,X)= $ $EAP_0(\jj,X)\oplus AP(\jj,X)$,\quad
$AAP(\jj,X)= C_0(\jj,X)$
$\oplus AP (\jj,X)$,

\noindent $ALAP_b(\jj,X)=$
 $ C_0(\jj,X)\oplus LAP_b (\jj,X)$\quad and $AAA(\jj,X)=$
 $ C_0(\jj,X)\oplus AA (\jj,X)$.

\smallskip

\noindent These are  the spaces consisting  respectively of the zero
function (when $\jj=\r$),  continuous functions vanishing at
infinity, almost periodic  ([1], [6], [24]), Levitan bounded  almost periodic [24], almost automorphic  functions [8], Eberlein (weakly)  almost periodic ([6, Definition 2.3.1]),
asymptotically almost
periodic functions  (when $\jj=\r_+$)([6, Definitions 2.2.1,
2.3.1, (2.3.2)]), asymptotically Levitan bounded almost
periodic functions and asymptotically almost
automorphic functions.

 For $\la\in \cc_+$  set

\smallskip

\qquad \qquad $f_{\la}(t) = \cases{ e^{-\la t},
\text{\,\, if\,} \,\, t \ge 0}
\\ { 0,\,\,\,\text {\,\,\,\,\,\,\,\,  if \,}\,\, t<
 0}\endcases$\qquad and $f_{-\la}= -\check {f_{\la}}$.

\smallskip

\noindent Then  $f_{\la}, \check {f_{\la}} \in L^1(\r)$ for all $\la\in \cc\setminus i\r$.
If $H \in L^{\infty} (\r,X)$  then

\noindent  $H* \check{f_{\la}}\in BUC(\r,X)$ for all $\la\in \cc\setminus i\r$.   We will  consider the property

\smallskip

  (3.7)\qquad $(\overline{F}* \check{f_{\la}})|\jj \in \A_{ub}$ for each $F\in \A$ and $\la\in \cc\setminus i\r$.

\proclaim {Proposition 3.1}  Let  $\A \st L^{\infty} (\jj,{X})$ be a closed subspace  satisfying (3.1) and (3.6). Assume  that $H\in W(\r,X)$ and $F = H| \jj$, where  $W(\r,X)$ and $V(\r)$ satisfy (2.4) and (2.8).

(i)    $sp_{\A,V} (H)= sp_{\A,V} (F)=  sp_{\A,V} (\overline{F})$.

(ii)  If  $F\in \A$,  then
  $(H*f)|\,\jj \in \A_{ub} $ for each $f\in V(\r)$.

(iii) If   $H\in L^{\infty} (\r,X)$, then
 $sp_{\A,\f} (H)= sp_{\A,\f} (F)=sp_{\A} (F)=$
  $sp_{\A} (H)$.

(iv) If  $ H\in \f'_{ar}(\r,X)$ and $0$ is an $(\A,\f)$-regular point
for $ H$,
  then there is $\delta >0$ and $\psi\in\f(\r)$ such that $\widehat{\psi}\in \h(\r)$, $\widehat{\psi} =1$ on $[-\delta,\delta]$ and $(H*\psi)|\,\jj \in \A_{ub}$.
\endproclaim

\demo{Proof} (i) If $\jj=\r$ there is nothing to prove so take $\jj=\r_+$. For $\va\in V(\r)$ we have $H * \va  =
\overline{F}*\va  + (H-\overline{F})*\va  $. By Lemma 2.3(ii), $((H-\overline{F}) *\va)|\r_+\in C_0 (\r_+,X)$, so by (3.2) it follows that $(H * \va) |\r_+ \in \A$ if and only if $(\overline{F} * \va) |\r_+ \in \A$.

(ii) Again we need only  consider the
   case $\jj=\r_+$.  For  $f\in V(\r)$ we have $(H*f) |\,\r_+=(\overline{F}*f) |\,\r_+ + \xi$ where $ \xi=((H-\overline{F})* f)|\, \r_+ \in C_0 (\r_+,X)$ as in part (i).  By (3.6), it follows that  $(\overline{F}*f) |\,\r_ + \in \A_{ub}$. Hence   $(H*f) |\,\r_+\in \A_{ub}$ by (3.2).

 (iii) By part (i) we have $sp_{\A,\f} (H)= sp_{\A,\f} (F)$ and   $sp_{\A,L^1 } (H)=sp_{\A,L^1} (F)$.
  Moreover,  $sp_{\A,L^1} (F)\st sp_{\A,\f} (F)$. For the reverse inclusion, let
   $\om_0\in \r$ be  $(\A,L^1)$-regular for $
F$.  So  there is $h_0\in L^1(\r)$ such
that $\widehat{h_0}(\om_0)\not =0$ and $(\overline{F}*h_0)|\,\jj \in \A$. Choose $\delta >0$ such that
 $\widehat{h_0}\not = 0$ on $
[\om_0-\delta,\om_0+\delta]$ and by   Lemma 2.1,
  $k_0\in L^1(\r)$ such that $\widehat{k_0}\cdot \widehat{h_0}=1
$ on $ [\om_0-\delta,\om_0+\delta]$. Let $\va\in \f (\r)$, $\widehat {\va}(\om_0)\not =0$ and
supp $\widehat {\va} \st [\om_0-\delta,\om_0+\delta]$.
 By (2.11)  we have $ \overline{F}*\va =
\overline{F}*(h_0*k_0*\va)=(\overline{F}*h_0)*(k_0*\va)$. So, $ (\overline{F}*\va)|\jj \in \A$ by part (ii) and therefore
 $\om_0$ is  $(\A,\f)$-regular for $F$.

(iv)  By (i), $0$ is   an $(\A,\f)$-regular point for $F$, so
 there is $\delta > 0$ and  $\va\in \f(\r)$ such that  $\widehat{\va}\not =0$ on $[-\delta,\delta]$ and $(\overline{F}*\va)|\,\jj \in \A$.  If $\jj=\r$, then $H=F=\overline{F}$ and so $H*\va \in \A$. If $\jj=\r_+$, then  $(H- \overline{F})*\va |\r_+ \in C_0(\r_+,X)$ by Lemma 2.3 (ii). So, $(H*\va)|\r_+ \in \A$ by   (3.2).
  By Lemma 2.1, there is $g\in \f (\r)$ such that $\widehat{g}\in \h(\r)$ and $\widehat{\va}\cdot \widehat{g} =1$ on $[-\delta,\delta]$. Obviously $\psi=\va*g\in \f(\r)$ and $\widehat{\psi}\in \h(\r)$.  Since $(H*\va)|\jj \in \A \st L^{\infty}(\jj,X)$, $(H*\psi)|\jj= ((H*\va) *g)|\jj\in \A_{ub} $ by (2.11) and part (ii).
$\P$
\enddemo

\proclaim{Proposition 3.2} Let  $\A\st L^{\infty}(\jj,X)$ be a closed subspace.

(i)  If $\A$ satisfies (3.6), then
  $\A$  is $BUC$-invariant.

(ii)   $\A$ satisfies (3.7) if and only if $\A$ satisfies (3.6).
\endproclaim
\demo{Proof}
(i) Take $F\in BUC(\r,X)$ with $F|\jj \in \A$ and take $t\in \r$. By Proposition 3.1(ii), for each $f\in L^1(\r)$ we have $(F_t*f)|\jj =(F*f_t)|\jj \in \A$. Using the approximate identity of Lemma 2.2 (i), we conclude that $F_t|\jj \in \A$.

(ii) Obviously, (3.6) implies (3.7). For the converse we begin by showing  that  $ E= span \,\{f_{\la}: \text{ Re\,\,}\la \not =0\}$ is a dense  subspace of  $L^1(\r)$.
    Indeed, if $E$ is not dense in $L^1(\r)$, then by the Hahn-Banach theorem there is
$0\not = \phi\in L^{\infty}(\r) = (L^1(\r))^*$ such that
$\Cal  {C}\phi(\la)= \int_0^{\infty} e^{-\la t}\phi (t)\, dt=0$
if  Re $\, \la >0$
and $\Cal  {C}\phi(\la)= -\int_0^{\infty} e^{\la t}\phi (-t)\, dt=0$ if Re $\, \la  < 0$. This means that the Carleman transform $\Cal  {C}\phi$ is zero on $\cc\setminus i\r$ (see (4.5) below) and implies $sp^{C} (\phi)=\emptyset$ and so $\phi =0$ by [28, Proposition 0.5 (ii)]. This is a contradiction  showing that $E$ is dense in $L^1(\r)$. Given (3.7) it follows
  that $(\overline{F}*f)|\, \jj \in \A$ for each $F\in \A$ and $f\in E$. Since $E$ is a dense subspace of $L^1(\r)$ and $\A$ is closed, (3.6) follows.
 $\P$
\enddemo

\proclaim {Remark 3.3} (a)
If  $\A\st BUC(\jj,X)$ satisfies (3.1) then using the properties of Bochner integration  (see [6, Lemma 1.2.1]) we find that  $\A$ satisfies (3.6).

 (ii) The converse of Proposition 3.2 (i) is false  in general.  The Banach space

 (3.8) \qquad $\A_{g} = g\cdot AP(\r,X)$ with
$g(t) =e^{it^2}$ for $t\in\r$

\noindent  satisfies (3.1) but does not satisfy (3.6).

  (iii) As  $\A_{g} \cap BUC(\r,X)=\{0\}$,  we conclude   that if $0\not = F\in BC(\r,X)$, then   by (3.5) and (1.1),

(3.9) \qquad  $sp_{\A_{g}} (F)= sp_0(F)= sp^B (F)\not = \emptyset$.

\noindent   In particular,
$sp_{\A_{g}} (F)\not = \emptyset$ for each $0\not =F\in \A_{g}$.
 \endproclaim

 \proclaim{Proposition 3.4}  Let  $\A \st L^{\infty} (\jj,{X})$ be a closed subspace  satisfying (3.6). Let   $W,V$ satisfy  (2.4),  (2.8) and $F, H \in W (\jj,X)$.

 (i)  If   $g\in V (\r)$ or $g\in L_c^{\infty} (\r)$,  then
 $sp_{\A,V} (\overline{F}*g) \st sp_{\A,V} (F) \,\cap $ supp $ \widehat
{g} $.

 (ii) $sp_{\A,V} (F)= \cup_{h > 0} sp_{\A,V} (M_h F)$.

(iii) If   $t\in \r$ and  $0\not = c\in \cc$, then
$sp_{\A,V} (c(\overline{F})_t) = sp_{\A,V} (F)$.

(iv)  $sp_{\A,V} (F+H) \st sp_{\A,V} (F) \cup sp_{\A,V} (H) $.

(v) If   $\gamma_{\la} \A \st \A$ for all $\la\in \r$,  then $sp_{\A,V} (\gamma_{\om}F)=\om +sp_{\A,V} (F)$ for all $\om\in \r$.
\endproclaim

\demo{Proof}
(i)  Assume
$\om\not \in sp_{\A,V} (F)$. Then there is $\va\in V (\r)$ with
$\widehat{\va} (\om)\not =0$ and $(\overline{F}*\va)|\jj \in \A$.   By (2.11),
 we have $(\overline{F}*g)*\va= (\overline{F}*\va)*g= \overline{F}*(\va*g)$. So, by
Proposition 3.1 (ii), we get  $((\overline{F}*\va)*g)|\,\jj\in \A$
proving $\om\not \in sp_{\A,V} (\overline{F}*g)$. On the other hand if $\om\not
\in $ supp $ \widehat {g}$, then there is $\va\in V (\r)$ with
$\widehat{\va} (\om)\not =0$ and  $\va *g=0$. So, $\om\not \in
sp_{\A,V} (\overline{F}*g)$. For the case $\A=\{0\}$ see also [28, Proposition
0.6 (i)].

(ii)  We note from (2.7) that $M_h F= (\overline{F}*s_h) |\jj$,  where $s_h \in L_c^{\infty} (\r)$  for each $h >0$. Hence $sp_{\A,V} (\overline{F}*s_h)\st sp_{\A,V} (\overline{F})$. By Proposition 3.1 (i), we have  $sp_{\A,V} (M_h F)=sp_{\A,V} (\overline{F}*s_h)$ and so
 $ \cup_{h > 0}sp_{\A,V} (M_h F)\st sp_{\A,V} (F)$. Now, let $\om\in sp_{\A,V} (F)$. There is $h >0$ such that $\widehat{s_h} (\om)\not = 0$. Assume that  $\om \not\in sp_{\A,V} (M_h F)= sp_{\A,V} (\overline{F}* s_h)$. There is  $\psi\in V(\r)$ such that $\widehat{\psi} (\om)\not = 0$ and $((\overline{F}*s_h)*\psi)|\jj \in \A$. By (2.11), $(\overline{F}*s_h)*\psi=\overline{F}*(s_h*\psi)$ so $(\overline{F}*(s_h*\psi))|\jj \in \A$. Since $s_h*\psi\in V(\r)$ and $\widehat {s_h*\psi} (\om)\not =0$, we conclude that $\om \not\in sp_{\A,V} (F)$, a contradiction which shows $\om \in sp_{\A,V} (M_h F)$.
 This  proves  $sp_{\A,V} (F)\st \cup_{h > 0}sp_{\A,V} (M_h F)$.

 The proofs of (iii), (iv) and (v) are similar to the case $\A=\{0\}$ ([28, Proposition 0.4]). \P
\enddemo

We recall (see [8, p. 118], [9, p. 1007],  [13], [31]) that a
function $F\in L^1_{loc} (\jj,X)$ is called $ergodic$  if there is
a constant $m(F)\in X$ such that

\qquad \qquad \qquad sup$_{t\in \jj}|| \frac{1}{T}\int_0^T
F(t+s)\,ds -m(F)||\to 0 $ as $T\to \infty$.

\noindent The limit $m(F)$ is called the $mean$ of $F$. The set of
all such ergodic functions will be denoted by $\E(\jj,X)$. We set $\E_{0}(\jj,X)=\{F\in \E(\jj,X): m(F)=0\} $,
$\E_b(\jj,X) =\E(\jj,X)\cap L^{\infty} (\jj,X)$,
$\E_{b,0}(\jj,X)=\{F\in \E_b(\jj,X): m(F)=0\} $, $\E_{ub}(\jj,X)
=\E(\jj,X)\cap BUC(\jj,X)$ and $\E_{u,0}(\jj,X)
=\E_{ub}(\jj,X)\cap \E_{b,0} (\jj,X)$.

 If $F\in L^1_{loc} (\jj,X)$ and
  $\gamma_{\om}F\in \E (\jj,X)$ for some $\om\in\r$, then

(3.10)\qquad $\gamma_{\om}M_h\, F\in \E (\jj,X)$ \,\, and \,\,
  $M_h\, \gamma_{\om} F\in \E_b (\jj,X)$\,\, for all $h >0$.

\noindent Moreover, if $F\in L^{\infty} (\jj,X)$ and
$\gamma_{\om}F\in \E_b (\jj,X)$ for some $\om\in\r$, then

(3.11) \qquad $\gamma_{\om}(\overline{F}*g)|\,\jj\in \E_{ub}(\jj,X)$ \,\,\, for all
$g\in L^1 (\r)$.

\noindent To prove (3.10), note that

 $M_T \gamma_{\om} M_h  F= \gamma_{\om}M_h
\gamma_{-\om}M_T \gamma_{\om} F$ and $M_T  M_h \gamma_{\om} F= M_h
M_T \gamma_{\om} F$.

\noindent It follows that $\gamma_{\om}M_h F$, $M_h
\gamma_{\om} F\in \E (\jj,X)$ for all $h >0$. By [9, (2.4)], $M_h
\gamma_{\om} F\in C_b (\jj,X)$ and so $M_h \gamma_{\om} F\in \E_b
(\jj,X)$. For (3.11) note that  if $F\in L^{\infty} (\jj,X)$, then
 $M_h F=(\overline{F}*s_h)|\jj$ (see  (2.7))  is  bounded and
uniformly continuous.  So, $\gamma_{\om} (\overline{F}*s_h)\in \E_{ub} (\jj,X)$
by(3.10).  A similar calculation gives $\gamma_{\om} (\overline{F}*\check {s_h})\in \E_{ub} (\jj,X)$. It follows that $ \gamma_{\om}(\overline{F}*g)|\jj \in \E_{ub} (\jj,X)$
for any step function $g$. Since step functions are dense in
$L^1(\r)$, (3.11) follows.

 Also, we note that

 (3.12)\qquad   $ \E_u(\jj,X): =UC(\jj,X)\cap
 \E(\jj,X)=\E_{ub}(\jj,X)$.

\noindent This follows by Lemma 2.2 (ii) using (3.10) (see also [9, Proposition 2.9]).

Next we recall the definition of the class of slowly
oscillating functions

\qquad $SO(\jj,X)= UC(\jj,X)+ L^1_{loc,0}(\jj,X)$,

\noindent where (see [18, Lemma 1.6], [5, Proposition 4.2.2] for the case $\jj=\r_+$)

\qquad $L^1_{loc,0}(\jj,X)=\{F\in L^1_{loc}(\jj,X): \lim_{|t|\to
\infty, t\in \jj}F (t)=0 \}$

\noindent
It follows that if $ F\in  L^1_{loc,0}(\jj,X)$ and  $\psi\in\f(\r)$, then

 (3.13)\qquad $  F\in \E_0(\jj,X)$ and  $ \overline{F}\in  L^1_{loc,0}(\r,X)$;

  (3.14)\qquad  $M_h F\in C_0(\jj,X)$  for all $h >0$ and $\overline{F}*\psi \in C_0(\r,X)$.

\noindent Also,  it is readily verified that for $ F\in SO(\jj,X)$

(3.15)\qquad  $\overline {F}*\psi, \,\,\, M_h \overline {F} \in UC(\r,X)$ for each $\psi\in \f(\r)$, $h > 0$.

\proclaim{Theorem 3.5}  Let $F\in \f'_{ar}(\r,X)$   and  $0\not\in sp_{0,\f} (F)$.

(i) If $PF \in UC(\r,X)$, then $PF \in BUC(\r,X)$ and
 $F\in \E_{0} (\r,X)$.

 (ii) If  $F\in L^{\infty}(\r,X)$ or  $F\in SO(\r,X)$,
 then $PF \in BUC(\r,X)$.
\endproclaim

\demo{Proof} (i) Choose $\va\in\f(\r) $ such that supp $\widehat{\va} \st [-\delta,\delta]$  and $\widehat{\va}=1$ on $[-\delta/2,\delta/2]$, where  $\delta > 0$ and  $ [-\delta,\delta] \cap sp_{0,\f} (F)=\emptyset$. Then $ sp_{0,\f} (F*\va) =\emptyset$  by Proposition 3.4 (i).  So $F*\va=0$ by [28, Proposition 0.5 (ii)]. Since $ (P F*\va)'= F*\va=0$, we  get $PF*\va=$ constant. By (2.12), we have $0\not \in sp_{0,\f} (PF- PF*\va)$. Hence $(PF- PF*\va)\in BUC(\r,X)$ by [12, Theorem 4.2] implying $PF\in BUC(\r,X)$. It is readily verified that $F\in \E_{0} (\r,X)$.

(ii) If $F\in L^{\infty}(\r,X)$ then clearly $PF \in UC(\r,X)$. So suppose that $F\in SO(\r,X)$ and  $h >0$. By (3.15) and Proposition 3.4 (ii) we have $M_h F \in UC(\r,X)$ and $0\in sp_{0,\f} (M_h F)$. Again $M_h F \in BUC(\r,X)$ by [12, Theorem 4.2]. As $\Delta_h PF = h M_hF$, one gets that $PF \in UC(\r,X)$ by [9, Proposition 1.4]. It follows that $PF \in BUC(\r,X)$ by part (i). \P
\enddemo

A slight modification of the  proof of Theorem 3.5 (i) gives the following  sharper result.
If $u\in UC(\r,X)$ and if $0\not \in $ supp $\widehat {u'}$, in the distributional sense,  then $u\in BUC(\r,X)$.
\smallskip

We are now ready to state and prove our main results.

\proclaim{Theorem 3.6} Let  $\A \st L^{\infty} (\jj,{X})$ be a closed subspace  satisfying (3.6) and $
\gamma_{\la}\A\st \tilde{\E}\in \{ \E(\jj,X), \E_{0}(\jj,X) \}$ for all $\la\in \r$. Let  $F\in \f'_{ar}(\jj,X)$  and $0\not\in sp_{\A,\f}(F)$.

(i)    $\overline{F}=H+G$, where  $H\in BUC(\r,X)$, $H|\,\jj\in \A_{ub}$ and $0\not\in sp_{0,\f} (G)$.

(ii) If  $PF\in UC(\jj,X)$ or  $F\in L^{\infty}(\jj,X)$ or  $F\in SO(\jj,X)$  then $F \in \tilde{\E}$. If also
    $F\in UC(\jj,X)$, then $ F\in \tilde{\E}\cap BUC(\jj,X)$.
\endproclaim

\demo{Proof} (i)  By Proposition 3.1 (iv)    there is  $\delta >0$ and $\psi\in
\f(\r)$ such that
 supp $\, \widehat {\psi}$ is compact,
 $\widehat {\psi} =1$ in a neighbourhood of $0$ and
$(\overline{F}*\psi)|\,\jj \in \A_{ub}\st \tilde{\E} $.
  Set $H= \overline{F}*\psi $ and $G=\overline{F}-H$. By (2.9), $H$ is continuous and  by Lemma 2.3, we conclude that $H\in BUC(\r,X)$. By  (2.12),  $0\not\in
sp_{0,\f} (G)$.

(ii) Since $PH \in UC(\r,X)$ it follows that $PG\in UC(\r,X)$ and by Theorem 3.5,
$G\in \E_0 (\r,X)$.  This and part (i) give $F = (H+G)|\jj \in \tilde{\E}$.
 The last assertion follows by  (3.12). \P
\enddemo

Consider the conditions

(3.16)\qquad $\A_0 (\jj,X)\in\{C_0 (\jj,X),EAP_0 (\jj,X)\}$, $\A(\jj,X) = \A_0(\jj,X) \oplus AP(\jj,X)$ and

(3.17)\qquad $\A_*(\jj,X) \in\{\A_0(\jj,X),\A(\jj,X)\}$.

\noindent In the notation of [6], $\A(\jj,X)$ is a $\Lambda$-class.

\proclaim{Theorem 3.7} Let $\A_0,\A, \A_*$  satisfy (3.16), (3.17). Assume that $F \in \f'_{ar}
(\jj,X)$,  $sp_{\A (\jj,X),\f}(F)$ is countable and  $\gamma_{-\om} F $
 $\in
\E(\jj,X)$ for all $\om\in sp_{\A(\jj,X),\f}(F)$.

(i) If $F \in UC(\jj,X)$, then $F\in \A(\jj,X)$. If also $sp_{\A_0(\jj,X),\f}(F)$ $=\emptyset$, then $F\in\A_0(\jj,X)$.

 (ii) If $F \in SO(\jj,X)$, then $F\in \A(\jj,X)+
L^1_{loc,0}(\jj,X)$.   If also $sp_{\A_0(\jj,X),\f}(F)$

\noindent $=\emptyset$, then $F\in
\A_0(\jj,X)+L^1_{loc,0}(\jj,X)$.

 (iii) If $F=H|\,\jj$ where $H \in L^{\infty}(\r,X)$ and if  $f\in L^1(\r)$, then $(H*f) |\jj\in \A(\jj,X)$.

(iv) If  $sp_{\A_*(\jj,X),\f}(F)=\emptyset$ and if $\psi\in \f(\r) $ with
$\widehat{\psi}\in \h(\r)$,
then $(\overline{F}*\psi)|\jj \in \A_*(\jj,X)$ and when $\jj=\r_+$, $(\overline{F}*\psi)|\r_- \in C_0(\r_-,X)$.

(v) If  $M_h F\in BC(\jj,X)$ for all $h> 0$ (for example if $F$ is ergodic) and if $\psi\in\f(\r) $, then $(\overline{F}*\psi)|\jj \in \A (\jj,X)$  and when $\jj=\r_+$, $(\overline{F}*\psi)|\r_- \in C_0(\r_-,X)$. If also $sp_{\A_0(\jj,X),\f}(F) = \emptyset$, then $(\overline{F}*\psi)|\jj \in \A_0(\jj,X)$.
\endproclaim
\demo{Proof} Assume that $F$ satisfies  the assumptions of one of the parts (i)-(iii).
We note that $\A(\jj,X)$ satisfies the assumptions of Theorem 3.6 with $\tilde{\E} =\E (\jj,X)$;
so, if $0\not\in
sp_{\A(\jj,X),\f}(F)$  then  $F$ is ergodic by Theorem 3.6.
If $0\in sp_{\A(\jj,X),\f}(F)$, then  $F$ is ergodic  by assumption.

(i) By (3.12) we get $F \in \E_{ub}(\jj,X)$.
   Let $\tilde{F}\in
BUC(\r,X)$ be an extension of $F$. By Proposition 3.1 (iii), $sp_{\A (\jj,X)}(\tilde{F})=sp_{\A (\jj,X)}({F})$ which is
countable.  By [6, Theorem 4.2.6],  $F =\tilde{F}|\jj \in
\A(\jj,X)$. If $sp_{\A_0(\jj,X),\f}(F)$
 $=\emptyset$, then  since $\A_0 (\jj,X)\st \E_{u,0} (\jj,X)$ and by Theorem 3.6(ii), $\gamma_{\la}F \in \E_{u,0}(\jj,X)$ for all $\la\in \r$. This implies $F\in
\A_0(\jj,X)$.

(ii)   Let  $F= u+\xi$, where $u\in UC(\jj,X)$, $\xi\in L^1_{loc,0}(\jj,X)$. We note that  $ sp_{\A(\jj,X),\f}(\xi)\st sp_{\A_0(\jj,X),\f}(\xi)=\emptyset$ by (3.14) and $\gamma_{\la}\xi \in \E_0(\jj,X)$ for all $\la\in \r$, by (3.13). Also, we have  $sp_{\A (\jj,X),\f}(M_h
F)$ is countable by Proposition 3.4 (ii).   By (3.10), we get $\gamma_{-\om}M_h F\in \E(\jj,X)$ for
all $\om \in sp_{\A(\jj,X),\f}(F)$. By Proposition 3.4(iv), $sp_{\A(\jj,X),\f}(M_h u)$ is countable  and  $\gamma_{-\om}M_h u\in \E(\jj,X)$ for
all $\om \in sp_{\A(\jj,X),\f}(F)$.
 So,  by part (i), we conclude that
$M_h u \in \A(\jj,X)$ for  all $h >0$. By Lemma 2.2 (ii), $ u =\lim _{h\to 0} M_h u\in \A(\jj,X) $. It follows that $F\in \A(\jj,X) + L^1_{loc,0}(\jj,X)$. If $sp_{\A_0(\jj,X),\f}(F)$
 $=\emptyset$, then  again $\gamma_{\la}F \in \E_{0}(\jj,X)$ for all $\la\in \r$, by Theorem 3.6 (ii). This implies that $F\in
\A_0 (\jj,X)+ L^1_{loc,0}(\jj,X)$.

(iii) Let $f\in L^1 (\r)$. Then
$\overline{F}*f \in BUC(\r,X)$. By Proposition 3.4(i), we deduce that
$sp_{\A(\jj,X),\f}(\overline{F}*f)$ is countable. By (3.11) we find that
$\gamma_{-\om}(\overline{F}*f)|\,\jj \in \E_{ub}(\jj,X)$ for all $\om\in
sp_{\A(\jj,X),\f}(\overline{F}*f)$. It follows  that $(\overline{F}*f)|\,\jj \in  \A(\jj,X))$, by part (i). By Lemma 2.3 (ii) and (3.2), we have $ ((H-\overline{F})*f)|\,\jj\in C_0 (\jj,X)\st \A_0 (\jj,X)$. Hence  $(H*f)|\,\jj \in  \A(\jj,X)$.

(iv)  Let $\om \in K=$ supp $\widehat{\psi}$. Since $\A_*(\jj,X)$ satisfies (3.1) and (3.6), by Proposition 3.1 (iv), there is
$f^{\om}\in \f(\r)$ such that $\widehat{f^{\om}}$ has compact
support, $\widehat{f^{\om}}=1$ on an open neighbourhood $V^{\om}$
of $\om$ and $(\overline{F} *f^{\om})|\,\jj \in \A_*(\jj,X)$. Take
$k^{\om}=f^{\om}*g^{\om}$, where $g^{\om} (t)=\overline{f^{\om}
(-t)}$. By  (2.11) and Proposition 3.1(ii),  we conclude that $(\overline{F} *k^{\om})|\,\jj \in \A_*(\jj,X)$.
Consider the open covering $\{V^{\om}: \om\in K\} $. By
compactness, there is a
 finite sub-covering $\{V^{\om_1},\cdots, V^{\om_n}\}$ of $K$.
 One has $k=\sum_{i=1}^n k^{\om_i}\in \f(\r)$, supp $\widehat{k}$ is compact,
 $\widehat{k}\ge 1$ on $K$ and $(\overline{F}*k) |\jj\in
\A_*(\jj,X)$. By Lemma 2.1, there is $h\in \f (\r)$ such that
$\widehat{h}\cdot \widehat{k}=1$ on $K$. Again by (3.6) and Proposition 3.1 (ii), it follows that
 $(\overline{F}*\psi) |\,\jj= ((\overline{F}*k)*h*\psi) |\jj\in \A_*(\jj,X)$. By Lemma 2.3 (i), if $\jj=\r_+$, then $(\overline{F}*\psi)|\r_-\in C_0((\r_-,X)$.

(v) Let $h >0$. By (3.10) and Proposition 3.4 (ii),
$\gamma_{-\om} M_h F\in \E_b(\jj,X)$ for all $h>0$ and $\om \in
sp_{\A (\jj,X),\f}(M_h F)\st sp_{\A (\jj,X),\f}(F)$. Therefore, by part (iii), $((\overline{F}*s_h) *g)|\jj \in
\A(\jj,X)$ for all $g\in L^1(\r)$. Take $\psi\in \f(\r)$. It
follows that $M_h(\overline{F} *\psi) |\jj=((\overline{F}*s_h) *\psi) |\jj\in \A(\jj,X)$ and also $(\Delta_h (\overline{F}*\psi))|\jj= (\overline{F}*\Delta_h
\psi)|\jj = (\overline{F}* hM_h \psi')|\jj =(\overline{F}*(s_h *\psi')) |\jj= ((\overline{F}*s_h)*\psi')) |\jj\in \A(\jj,X) \st BUC(\jj,X)$.  By
[9, Proposition 1.4], one gets $(\overline{F}*\psi)|\,\jj$ is uniformly continuous.
This implies $(\overline{F}*\psi)|\jj= \lim_{h\searrow  0} M_h(\overline{F}*\psi)|\jj\in \A(\jj,X)$,  by Lemma 2.2 (ii). If $\jj=\r_+$ we proceed as in (iv).  If   $sp_{\A_0(\jj,X),\f}(F) = \emptyset$, the result follows by (i). \P
\enddemo

In the following we demonstrate again how one can extend the results proved for the case $F\in BUC(\r,X)$ to the case $F\in \f'_{ar} (\r,X)$.

\proclaim{Theorem 3.8} Assume that $F \in \f'_{ar}
(\r,X)$,  $sp_{AP ,\f}(F)$ is countable and  $\gamma_{-\om} F $
 $\in
\E(\r,X)$ for all $\om\in sp_{AP,\f}(F)$.

 (i) If $F \in UC(\r,X)$ , then $F\in AP(\r,X)$.

(ii) If $F \in SO(\r,X)$, then  $F=u$  almost everywhere for some $u\in AP(\r,X)$.

 (iii) If $F \in L^{\infty}(\r,X)$ and if  $\psi\in L^1(\r)$, then $F*\psi \in AP(\r,X)$.

(iv) If  $sp_{AP,\f}(F)=\emptyset$ and if $\psi\in \f(\r) $ with
$\widehat{\psi}\in \h(\r)$,
then $F*\psi \in AP (\r,X)$.

(v) If $M_h F\in BC(\r,X)$ for all $h> 0$,
then $F*\psi \in AP(\r,X)$ for each  $\psi\in \f(\r)$.

\endproclaim

\demo{Proof} Since $AP(\r,X)\st AAP(\r,X)$, we conclude that $sp_{AAP ,\f}(F)\st sp_{AP ,\f}(F)$.

 (i) This follows by Theorem 3.7 (i) and [6, Theorem 4.2.6].

(ii)  Let $ h >0$.  By (3.15)  $M_h F$ is uniformly continuous. By (i), $M_h F\in AP(\r,X)$.
  So $sp_{AP,\f}(F) = \emptyset$ by (3.6) and Proposition 3.4(ii).
 By Theorem 3.7 (ii)  with $\A_0= C_0 (\r,X)$, one has $F=u+\xi$  where $ u\in AP(\r,X)$ and $\xi\in L^1_{loc,0} (\r,X)$.
 By Proposition 3.4 (iv), we get $sp_{AP,\f} (\xi)=sp_{AP,\f} (F-u) \st sp_{AP,\f}(F)\cup sp_{AP,\f} (-u)=\emptyset$.
By (3.14), if $\va\in \f(\r)$ then
   $\xi *\va \in AP(\r,X)$  if and only if  $\xi *\va= 0$. It follows that $sp_{0,\f} (\xi)=sp_{AP,\f} (\xi)= \emptyset $.   Hence $\xi =0$ (almost everywhere) by [28, Proposition 0.5 (ii)].

 (iii)-(v): By  Theorem 3.7 (iii)-(v), $F*\psi \in AAP (\r,X)$. Hence   $F*\psi \in AP (\r,X)$ by Proposition 3.4 (i) and part (i) above.
 \P
\enddemo

\proclaim{Corollary 3.9} Assume that $F \in \f'_{ar}
(\jj,X)$,  $\psi\in \f(\r) $ with  $\widehat{\psi}\in \h(\r) $  and $sp_{C_0(\jj,X),\f}(F)\,\, \cap $ supp $\widehat{\psi} =\emptyset$.  Then $\overline{F}*\psi \in C_0(\r,X)$.
\endproclaim
\demo{Proof} Choose $\rho\in \f(\r)$ such that $\widehat {\rho}\in\h(\r)$ and  $\widehat {\rho}=1$ on an open  neighbourhood of supp $\widehat {\psi}$. Then
$(\overline{F}*\psi)*\rho= \overline{F}*(\psi*\rho)=\overline{F}*\psi$. By Proposition 3.4 (i), $sp_{C_0(\jj,X),\f}(\overline{F}*\psi)=\emptyset$, so $(\overline{F}*\psi)*\rho\in C_0 (\r,X)$ by Theorem 3.7 (iv). \P
\enddemo
 \proclaim{Corollary  3.10} Assume $ F\in L^1_{loc}(\jj,X)$ and  $sp_{C_0 (\jj,X),\h} (F)=\emptyset$. If $(\overline{F}*\psi)|\jj$ is  uniformly continuous   for some
$\psi\in \h (\r)$,  then
 $\overline{F}*\psi\in C_0 (\r,X)$.
\endproclaim
\demo{Proof}   By $\h(\r)\st \f(\r)$ and Proposition 3.4 (i),

  $ sp_{C_0 (\r,X),\f}
(\overline{F}*\psi) \st sp_{C_0 (\r,X),\h}
(\overline{F}*\psi)  \st sp_{C_0 (\r,X),\h}
({F})=\emptyset$.

\noindent So,  $sp_{C_0 (\r,X),\f}
(\overline{F}*\psi)=\emptyset$. The result follows from Theorem 3.7 (i).
 \P
\enddemo

The following example shows that the assumption of uniform continuity is essential in Corollary 3.10.
\proclaim{Example 3.11} If $F (t)=e^t$ for  $t\in \r $, then  $sp_{C_0 (\jj,X),\h} (F)=\emptyset$   but $(F*\psi)|\jj$ is unbounded for  each  $\psi\in \h (\r)$ with  $
\int_{-\infty}^{\infty} e^{-s}\psi(s)\, ds \not = 0$.
\endproclaim
\demo{Proof}
 For any $\om\in \r$, choose $a >0$ such that
 $\cos\,
\om\,t$ does not change sign on $[0,a]$. Take $\va\in \h(\r)$ such
that $\va
>0 $ on $ (0,a)$ and supp $\va = [0,a]$. Let $f(t) =\va(t)$ for
$t\ge 0$, $f(t) = - e^{2t}\va(-t)$ for $t < 0$.
  It follows that $f\in \h(\r)$, $F*f =0$ and $\widehat{f}(\om)\not =0$. This means $sp_{C_0 (\jj,X),\h} (F)=\emptyset$.
  Moreover, for $\psi\in \h(\r)$  we have $F*\psi (t)= c e^t$, where $c=
\int_{-\infty}^{\infty} e^{-s}\psi(s)\, ds$.  So,
$(F*\psi)|\jj$ is unbounded if $c\not = 0$. \P
\enddemo

In the following example we calculate  reduced spectra
 of some functions  whose Fourier transforms  may not be regular distributions.

 \proclaim{Example 3.12}  (i) If $ F\in L^p(\jj,X)$ for some $1\le p
< \infty$, then $M_h F\in C_0(\jj,X)$ for all $h>0$ and
$sp_{C_0 (\r,X),V} (F)=\emptyset$ for any  $V\in \{\h (\r),\f(\r)\}$.

 (ii) Let $F\in
 \E_{ub}(\jj,X)$ and either $F'\in
L^p(\jj,X)$ for some $1\le p < \infty$ or  more generally  $F' \in L^1_{loc}(\jj,X)$ with $M_h
F'\in C_0(\jj,X)$ for all $h >0$. Then $F\in X\oplus C_0 (\jj,X)$ and  $sp_{C_0 (\jj,X)}(F)\st
\{0\}$.
\endproclaim
\demo{Proof} (i)   By H\"{o}lder's inequality,
 $||M_h F (t)||=(1/h )||\int_0^h F (t+s)\, ds||\le
h^{-1/p}$

 \noindent $(\int_0^h ||F(t+s)||^p\,ds)^{1/p})$, so
   $ M_h F \in C_0 (\jj,X)$
for all  $h >0$.  By (2.7) and Lemma 2.3 (i), we get $  \overline{F}*s_h \in C_0 (\r,X)$
for all  $h >0$. So, $sp_{C_0 (\r,X),V} (\overline{F}*s_h)=sp_{C_0 (\jj,X),V} (M_h F)=\emptyset$ for all $h
>0$. Hence
 $sp_{C_0 (\r,X),V} (F)=\emptyset$ by Proposition 3.4(ii).

(ii) By part (i)  we have $h M_hF' (\cdot)= F (\cdot+h)-F(\cdot)\, \in C_0
(\jj,X)$ for all $h>0$.  Let $\tilde {F}\in BUC(\r,X)$ be given by  $\tilde {F} =F$  on $\jj$ and $\tilde {F} (t) =F (0)$ on $\r\setminus \jj$.  It follows that $\Delta_s \tilde {F} \in C_0
(\r,X)$ for all $s\in\r$.  By [6, Theorem
4.2.2, Corollary 4.2.3], we conclude that $F= \tilde{F}|\jj\in X\oplus C_0 (\jj,X)$. This implies  $sp_{C_0 (\jj,X)}(F)\st
\{0\}$.  $\P$
\enddemo

 The following result shows that the ergodicity condition in Theorem 3.7 parts (i), (ii), (v) is necessary.
\proclaim{Example 3.13}  Let $F \in  C_0
 (\jj,X)$   with  $PF$  unbounded. Then  $sp_{C_0(\jj,X),\f} (PF)=\{0\}$ and
 $ P\overline{F}*\psi|\jj \not \in C_0(\jj,X)$ for each $\psi \in \f(\r)$ with $\widehat {\psi} (0)\not =0$.
\endproclaim
\demo{Proof} Note that $ P\overline{F}= \overline{PF}\in UC(\r,X)$.
Set $\va=\psi'$ where $\psi\in \f(\r)$. Then $P\overline{F}*\va|\jj=\overline{F}*\psi|\jj \in C_0(\jj,X)$. This shows that $sp_{C_0(\jj,X),\f} (PF)\st \{0\}$. If   $sp_{C_0(\jj,X),\f} (PF)=\emptyset$ we conclude that ${PF}\in C_0(\jj,X)$ by Theorem 3.7(i).  But $PF$ is  assumed to be unbounded, so
 $sp_{C_0(\jj,X),\f} (PF)=\{0\}$.
 Now, let $\psi \in \f(\r)$ with $\widehat {\psi} (0)\not =0$. If  $P\overline{F}*\psi|\jj  \in C_0(\r,X)$, then $0\not \in sp_{C_0(\jj,X),\f} (PF)$, a contradiction which proves $P\overline{F}*\psi|\jj \not \in C_0(\r,X)$. \P
\enddemo

 \head{\S 4. Properties of the weak Laplace spectra}\endhead

In this section we  establish some new properties of the  Laplace and weak
Laplace spectrum for  regular tempered distributions and show that
they are similar to those of the Carleman  spectrum (see [28,
Proposition 0.6]). We use the functions $e_a$ for $a \ge 0$ defined on $\r$ or $\r_+$ by $e_a (t)= e^{-at}$.

If $F\in  \f'_{ar}(\r_+,X)$, then
  $ e_a F\in L^1(\r_+,X)$ for all
$a> 0$
 and so the $Laplace$ $ transform$ $\Cal {L}F$  may be defined by

\smallskip

(4.1) \qquad $\Cal {L}F(\la)=\int_0^{\infty}\, e^{-\la \,t}F (t)\,
dt$ \,\,\,\, for\,\,\,\, $\la \in \cc_+ $.

\smallskip

\noindent For a function $F\in  \f'_{ar}(\r,X)$  the $Carleman$ $transform$  $\Cal {C}F$ is defined  by

\smallskip

(4.2) \qquad   $\Cal {C} {F} (\la)=
 \cases { \Cal{L^+}F (\la)= \int_0^{\infty}\, e^{-\la \,t}F (t)\,
dt \qquad \,{\text{for\,\,} \la
 \in \cc_+}}\\
  { \Cal{L^-}F(\la) = - \int _0^{\infty}\, e^{\la\, t} F(-t)\, dt\qquad {\text{\,
  for\,\,}
 \la \in \cc_-}.}\endcases$

\smallskip

If $F\in L^1(\r_+,X)$, then  $\Cal {L} F$ has a continuous
extension  to $\cc_+\cup i\r$  given also by the integral in
(4.1). By the Riemann-Lebesgue lemma $\widehat{\overline{F}} =\Cal
{L}F(i\cdot)\in C_0 (\r,X)$.

 If $F\in  \f'_{ar}(\r_+,X)$, then
 $\widehat {\overline{F}}\in  \f'(\r,X)$ and
$\Cal{L} F(a+i\cdot)=\widehat{e_a\overline{F}}  \in \f'_{ar}(\r,X)$ for all $a >0$.
Moreover, for $\va\in \f(\r)$,

\smallskip

(4.3)\qquad  $<\Cal {L} F (a+i\cdot),\va>\,\,= <\widehat{e_a\overline{F}},\va> =$
 $ <{e_a\overline{F}},\widehat {\va}>\,\,\to
<\overline{F},\widehat {\va}>\,\,= <\widehat {\overline{F}},\va>$,

\noindent where the limit
exists as $a\searrow 0$ by the Lebesgue  convergence
theorem. This means that $\lim_{a\searrow 0}\Cal{L} F (a+i\cdot)= \widehat{\overline{F}} $ with respect to the weak dual topology on $ \f'(\r,X)$.

For a  holomorphic function  $\zeta: \Sigma\to X$, where
$\Sigma=\cc_+$ or $\Sigma=\cc\setminus i\,\r$,  the point
$i\,\om\in i\,\r$ is called a $regular$ $point$ for $\zeta$ or
$\zeta$ is called $holomorphic$ at $i\,\om$, if $\zeta$ has an
extension $\tilde{\zeta}$ which is holomorphic in a
neighbourhood $V\st \cc$ of $i\,\om$.

 Points $i\,\om$ which are $not$ $regular$ $points$ are
called $singular$ $points$.

  The  $Laplace$ $spectrum$ of a function    $F\in  \f'_{ar}(\r_+,X)$   is defined by

\smallskip

 (4.4) \qquad $sp^{\Cal{L}}
(F) =\{\om\in\r: i\, \om$ is a singular point for $\Cal {L}F \}$.

\smallskip

\noindent  The $Carleman \,\, spectrum$ of a function    $F\in  \f'_{ar}(\r,X)$  is defined (see [5, (4.26)]) by

\smallskip

(4.5) \qquad $sp^{\Cal{C}} (F) =\{\om\in\r: i\, \om$ is a
singular point for $\Cal {C} {F}\}$.
\smallskip

\noindent The Laplace spectrum is also called  the half-line
spectrum ([5, p. 275]).

 Note that if $ \tilde{\Cal{L}} \gamma_{-\om} F $ and $
\tilde{\Cal{C}} \gamma_{-\om} F$  are  holomorphic extensions of
$ \Cal{L} \gamma_{-\om} F$ and $\Cal{C} \gamma_{-\om} F$
respectively, which are holomorphic in a neighbourhood of $0$, then

\smallskip

(4.6)\qquad $\lim _{\la\to 0}\Cal{L} \gamma_{-\om }F (\la)=
\tilde{\Cal{L}} \gamma_{-\om} F (0)$  if $\om \not\in
sp^{\Cal{L}} (F) $, and

\qquad\qquad\, $\lim _{\la\to 0}\Cal{C}
\gamma_{-\om }F(\la)= \tilde{\Cal{C}} \gamma_{-\om} F (0)$ if
$\om \not\in sp^{\Cal{C}} (F) $.

\smallskip

\noindent If  $F\in L^{\infty} (\r_+,X)$    and $
sp^{\Cal{L}}(F)=\emptyset$, then by Zagier's result [33, Analytic
Theorem] we conclude that $\widehat {F}(\om)=\int_0^{\infty}
e^{-i\om \,t} F (t)\,dt$ exists as an improper  integral (and by (4.6)
equals $\tilde{\Cal{L}} \gamma_{-\om} F (0)$) for each $\om \in
\r$. Zagier's analytic theorem does not hold  for unbounded
functions. Indeed, the Laplace spectrum of $F (t) = t e^{it^2}$ is
empty  (see Example 4.5 below) and  it can be verified that
$\int_0^{\infty}  e^{-i\om \,t} F (t)\,dt$ does not exist as an
improper Riemann integral for any $\om\in \r$.

For a  holomorphic function  $\zeta: \cc_+\to X$,  the point
$i\,\om\in i\,\r$ is called a  $weakly$ $regular$ $point$ for
$\zeta$ if there exist $\e
> 0$ and  $h\in L^1 (\om-\e, \om+\e)$ such that

\smallskip

(4.7)\qquad  $\lim_{a\searrow 0} \int _{-\infty}^{\infty}\zeta
(a+i\,s)\va (s)\, ds= \int _{\,\om-\e}^{\,\om+\e}h (s)\va (s)\,
ds$

\qquad  \qquad \,\,\, for all $\va\in \h(\r) $ with supp$
\,\va \st (\om-\e, \om+\e)$.

\smallskip

\noindent See [4, p. 474] for the particular case $h\in C(\om-\e, \om+\e)$. The  points $i\,\om$ which are not weakly regular points are
called $weakly$ $singular$ $points$.

 The $weak\,\, Laplace\,\, spectrum$ of $F\in
   \f'_{ar}(\r_+,X)$ is defined  (see [5, p. 324]) by

\smallskip

(4.8) \qquad $sp^{w\Cal{L}} (F) =\{\om\in\r: i\, \om$ is not a
weakly regular point for $\Cal {L} {F}\}$.

\noindent For  $F\in
   \f'_{ar}(\r,X)$, we define  $sp^{w\Cal{L}} (F):= sp^{w\Cal{L}} (F|\r_+)$. It follows readily that if $F\in
   \f'_{ar}(\r,X)$ then

\smallskip

(4.9) \qquad $sp^{w\Cal{L}} (F)\st sp^{\Cal{L}} (F)\st
sp^{\Cal{C}} (F)$; and,
 if $F\in L^1(\r_+,X)$,
  $sp^{w\Cal{L}} (F)=\emptyset$.

\smallskip

In the following $sp^{*}$ denotes $sp^{\Cal{L}}$ or $sp^{w\Cal{L}}$  or
$sp^{\Cal{C}}$. Note that $sp^{*} (F)$ is closed for any $ F\in \f'_{ar} (\r,X)$.

\proclaim{Proposition 4.1}  If $F, G
\in \f'_{ar}(\r,X)$, then

(i) $sp^{*} (F) =sp^{*} (F_s)= sp^{*} (cF)$ for each $s\in \r, 0\not = c\in \cc$.

(ii) $sp^{*} (F)=\cup_{h>0} sp^{*} (M_h F)$.

(iii)  $sp^{*} (\gamma_{\om}F)= \om + sp^{*} (F)$.

(iv) $sp^{*} (F+G) \st sp^{*} (F)\, \cup\, sp^{*} (G)$.
\endproclaim

\demo{Proof}  (i) A simple calculation shows that  for $\la\in
\cc^{\pm}$

\smallskip

(4.10)\qquad $\Cal {L^{\pm}}F_s (\la) = e^{\la \,s}\Cal {L^{\pm}}F
(\la)- e^{\la \,s} \int_0^{s}\, e^{-\la \,t }F (t)\,dt$.

\smallskip

\noindent Note that the second term on the right of (4.10) is entire in $\la$ for each $s\in \r$. It follows that $\Cal {L^+}F$  (respectively $\Cal {C}F$) is
holomorphic at $i\,\om$ if and only if $\Cal {L}^+ F_s $
(respectively $\Cal {C}F_s$) is holomorphic at $i\,\om$. This
proves (i)   for $sp^{\Cal{L}}$ and $sp^{\Cal{C}}$.  Now, assume $i\,\om$
is a weakly  regular point for $\Cal{L}F$. So there exists  $\e >0$
and $h \in L^1 (\om-\e,\om+\e)$ satisfying $\lim_{a\searrow 0}
\int _{-\infty}^{\infty}\Cal{L^+}F (a+i\,\eta)\va (\eta)\, d\eta=
\int _{\om-\e}^{\om+\e}h (\eta)\va (\eta)\, d\eta$
 for all $\va\in \h(\r) $ with supp$ \,\va \st
(\om-\e, \om+\e)$. Then by [30, Theorem 6.18, p. 146] (valid also for $X$-valued distributions),
$\lim_{a\searrow 0} \int _{-\infty}^{\infty}\Cal{L}F
(a+i\,\eta) e^{(a+i\,\eta) s}\va (\eta) d\eta= \int
_{\om-\e}^{\om+\e}h (\eta)e^{i\eta \,s} \,\va (\eta) d\eta$
 for all $\va\in \h(\r) $ with supp$ \,\va \st \, (\om-\e, \om+\e)$.
 It follows that $i\,\om$ is a weakly regular point for
$\Cal{L}F_s$.

(ii) Another  calculation shows that for  $\la \in \cc^{\pm}$

\smallskip

(4.11)\qquad $\Cal {L^{\pm}}M_h F (\la) =$
  $g(\la\, h) \Cal {L^{\pm}}F (\la) - (1/h)\int_0^h (e^{\la \, v} \int_0 ^v e^{-\la \, t}F(t)
 dt)\,dv$,

\smallskip

\noindent where  $g$ is the entire function given by
$g(\la)=\frac{e^{\la}-1}{\la}$ for $\la\not = 0$. Let $i\,\om\in
i\, \r$ be a regular point for $\Cal {L^+}F$ and let
$\tilde {\Cal {L^+}}F: V \to X$ be a holomorphic extension of $\Cal
{L^+}F$ to a neighbourhood $V \st \cc$ of $i\,\om$. Then
$\tilde{\Cal{L^+}}M_h F (\la) =$
  $g(\la\,h)\tilde  {\Cal {L^+}}F (\la) -
  (1/h)\int_0^h (e^{\la \, v} \int_0 ^v e^{-\la \, t}F(t)\,
 dt) \,dv$, $\la \in V$,  is a holomorphic extension of
$\Cal {L^+}M_h F$. So $i\,\om$ is a regular point for $\Cal {L^+}
M_h F$.
 Conversely suppose  $i\om\in i\r$  is a regular point of $\Cal{L} M_h F$ for each $h> 0$. Choose $h_0 >0$ such
that $g(i\om\, h_0)\not =0$. Then $i\,\om$ is
a regular point for $\Cal {L^+} F$. This proves (ii) for
$sp^{\Cal{L}}$. The case $sp^{\Cal{C}}$ follows similarly noting
that (4.11) implies $\Cal {C} M_h F (\la) =$
  $g(\la\, h) \Cal {C}F (\la) - (1/h)\int_0^h (e^{\la \, v} \int_0 ^v e^{-\la \, t}F(t)
 dt)\,dv$. The proof for  $sp^{w\Cal{L}}$ is
similar to the one in part (i).

(iii) This follows easily from the definitions noting that $\Cal
{L}^+(\gamma_{\om}F) (\la)= \Cal {L}^+ F (\la -i\om)$  and $\Cal
{C}(\gamma_{\om}F) (\la)= \Cal {C}F (\la -i\om)$.

(iv) This  follows directly from the definition.
 \P
\enddemo

The following result was obtained in  [10, (3.12)] in the case  $F \in L^{\infty}(\r_+,X)$ since then $sp_{C_0(\r_+,X),\f}(F)=sp_{C_0(\r_+,X)}(F)$ (see also [17, Lemma 1.16] for $\A= C_0(\r_+,X)$).

\proclaim{Proposition 4.2}  If $F \in \f'_{ar}(\r_+,X)$ and $\A\st L^{\infty}(\r_+,X)$ satisfies (3.1) then $ sp_{\A,\f}(F) \st sp_{C_0(\r_+,X),\f}(F) \st sp^{w\Cal{L}}(F)$.
\endproclaim
 \demo{Proof}  By  (3.2),  $C_0(\r_+,X)\st \A$ and so $ sp_{\A,\f}(F) \st sp_{C_0(\r_+,X),\f}(F) $.  Let $\om \not \in  sp^{w\Cal{L}}(F) $. Choose $\e >0$ and
$\va\in\f(\r)$ such that $sp^{w\Cal{L}}(F)\cap
[\om-\e,\om+\e]=\emptyset$, $\widehat{\va}(\om)=1$ and supp $\widehat{\va}\st
[\om-\e,\om+\e]$. By [18, Proposition 1.3], $\overline{F}*\va\in C_0(\r,X)$ and so $\om\not\in sp_{C_0(\r_+,X),\f}(F)$. \P
\enddemo

\proclaim{Theorem 4.3} Let   $F \in \f'_{ar}
(\r_+,X)$.

(i) If $0\not \in sp^{\Cal {L}}(F)$ and $PF \in UC(\r_+,X)$, then $PF \in BUC(\r_+,X)$ and  $F \in \E_0(\r_+,X)$. If also   $F\in L^{\infty}(\r_+,X)$ or  $F\in SO(\r_+,X)$,
 then $PF \in BUC(\r_+,X)$.

(ii) If  $M_h F\in BC(\jj,X)$ for all $h> 0$,   if  $sp^{w\Cal {L}}(F)$ is   countable and  $\gamma_{-\om} F $
 $\in
\E(\jj,X)$ for all $\om\in sp^{w\Cal {L}}(F)$,   and if $\psi\in\f(\r) $,
 then  $(\overline{F}*\psi) |\r_+\in
AAP(\r_+,X)$  for each $\psi\in \f(\r)$.
If also $ sp^{w\Cal {L}}(F)= \emptyset$,
then $\overline{F}*\psi \in C_0(\r,X)$.

\endproclaim
 \demo{Proof} (i) Let $h >0$. We have $M_h F= (1/h) \Delta_h PF \in BUC(\r_+,X)$ and $0\not \in sp^{\Cal {L}}(M_h F)$ by Proposition 4.1 (ii). By [5, Corollary 4.4.4], $ P M_h F \in BUC(\r_+,X)$. As  $  M_h P F = P M_h F$, by Lemma 2.2 (ii) we conclude that $PF \in BUC(\r_+,X)$ and hence $F \in \E_0(\r_+,X)$.

 If  $F\in L^{\infty}(\r_+,X)$, then clearly  $PF \in UC(\r_+,X)$. So, assume that   $F\in SO(\r_+,X)$. Then
  $M_h F \in UC(\r_+,X)$ by (3.15) and $0\not \in sp^{\Cal {L}}(M_h F)$ by Proposition 4.1 (ii). It follows that
   $M_h F \in BUC(\r_+,X)$ by [18, Proposition 1.3, Remark 1.4].  This implies $PF \in UC(\r_+,X)$ by [9, Proposition 1.4].

 (ii) By Proposition 4.2, we have $sp_{AAP,\f} (F) \st sp_{C_0,\f} (F) \st sp^{w\Cal {L}}(F)$ and the result follows by Theorem 3.7 (v) with  $\A =AAP (\r_+,X)$ and $\A_0 =C_0 (\r_+,X)$. \P
\enddemo

\noindent {\bf{Remark 4.4}}. (i)  In the case $\jj=\r_+$,  Theorem 3.6 and Theorem 3.7 parts (i)- (iv)  remain valid if we replace $sp_{\A,\f} (F)$ and $sp_{C_0(\r_+,X),\f} (F)$ by $sp^{\Cal{L}} (F)$ or $sp^{w\Cal{L}} (F)$.
 Indeed,  note that Theorem 3.6  holds for $\A= C_0 (\r_+,X)$. By   Proposition 4.2 and (4.9), we have $sp_{C_0(\r_+,X),\f} (F)\st sp^{w\Cal{L}} (F) \st sp^{\Cal{L}} (F)$. So, Theorem 3.6, Theorem 3.7 parts (i)- (iv) and Corollary 3.9
  strengthen the results of Chill  [17, Lemma 1.16], [18, Proposition 1.3, Theorem 1.5,  Corollary 1.7].
  Theorem 3.7 (v) and Theorem 4.3  seem to be new for any spectrum.

(ii) If $F$ in Theorem 3.6  is
not bounded or slowly oscillating, then $F$ is not necessarily ergodic.
For example, if $g (t)= e^{it^2}$ and $F= g^{(n)}$ for some $n\in \N$, then by
Example 4.5 below  and (4.9), we find $sp^{w\Cal{L}} (F)=\emptyset$. By
Proposition 4.2, we get $sp_{C_0 (\r_+,\cc),\f} (F )=\emptyset$ but
$F|\r_+ $ is neither bounded nor  ergodic  when $n\ge 2$. If $n=1$,
$F$ is ergodic but  not bounded.

 (iii) In view of Proposition 4.2 and (4.9) several tauberian theorems by Ingham  [22] ([5, Theorem 4.9.5]) and their
generalizations in [2], [3], [5, Theorem 4.7.7, Corollary 4.7.10, Theorem 4.9.7, Lemma 4.10.2], [14], [15], [17, Lemma 1.16, p. 25], [18] are consequences of
Theorem 3.6 and Theorem 3.7. Our
proofs are simpler and  different. Replacing Laplace and weak Laplace spectra by reduced spectra we are able to strengthen and unify these previous results.

(iv) If $F\in C_0 (\jj,X)$ satisfies $F= \frac{1}{|t|}$ for $t\in \jj$, $|t| > 1$, then $PF \in UC(\jj,X)$ but $PF$ is not bounded. This means that Theorem 3.5(i) and Theorem 4.3(i) are not valid if we replace $sp_{0,\f} (F)$ or $sp^{\Cal{L}} (F)$ by $sp_{C_0(\jj,X),\f} (F)$.

\smallskip

In the following  we use our results to calculate some  Laplace spectra.

\proclaim{Example 4.5} Take $g (t)=e^{it^2}$ for  $t\in \r$. Then  $sp^{\Cal {C}}(g) =\r$ and
$sp^{\Cal {L}}(g)= sp^{\Cal{L}}(g^{(n)}) =\emptyset$ for any
$n\in \N$. Moreover, $M_h g \in C_0
(\r,\cc)$ and $sp^{\Cal{L}}(M_h g) =\emptyset$ for all  $h>0$.
\endproclaim
\demo{Proof} By Proposition 4.1 (i), (iii), it is readily verified
 that $sp^{\Cal{L}}(g)=sp^{\Cal{L}}(g_a)=2a+ sp^{\Cal
 {L}}(g)$ for each  $a\in \r$. This implies that either $sp^{\Cal{L}}(g)
 =\emptyset$ or $sp^{\Cal{L}}(g) =\r$. Similarly either $sp^{\Cal {C}}(g)
 =\emptyset$ or $sp^{\Cal {C}}(g) =\r$.  But $g\not =0$ and so by [28, Proposition 0.5 (ii)], $sp^{\Cal {C}}(g)=\r$.

  Next note that $ y (\la)= \Cal
 {L}^+g (\la)$ is a solution of the differential equation $y'(\la) + (\la/2i)y (\la)= 1/2i$ for $\la\in \cc_+$. Solving the equation we find $y (\la)= e^{-\la ^2/4i} (c+ (1/2i) \int_0^{\la} e^{z^2/4i}\,dz$ for some choice of $c\in \cc$. As this last function is entire we conclude that $ sp^{\Cal{L}}(g)=\emptyset $.
 Since  $\int_0^{\infty} e^{i\,t^2}\, dt$ converges as
an improper Riemann integral and $M_h g(t)= (1/h)[Pg
(t+h)-Pg (t)]$ it follows that  $M_h g\in C_0 (\r,\cc)$ for each $h>0$. Moreover, by Proposition 4.1,  $sp^{\Cal{L}}(M_h g)\st sp^{\Cal{L}}(g)
=\emptyset$ and $sp^{\Cal{L}}(g')= \cup_{h >0} sp^{\Cal{L}}(M_h g')=\cup_{h >0} sp^{\Cal{L}}((1/h)(g_h-g))=\emptyset$. It follows that $sp^{\Cal{L}}(g^{(n)}) =\emptyset$ for any
$n\in \N$.
 \P
\enddemo

Finally, we demonstrate  that our results can be used to deduce spectral criteria for bounded solutions of evolution equations  of the form

\smallskip

 (4.12)\qquad  $\frac{d u(t)}{dt}= A u(t) +  \phi (t) $, $u(0) \in X$, $t\in  {\jj}$,

\smallskip

\noindent  where
 $A$ is a closed linear operator   on $X$ and $\phi\in L^{\infty} (\jj, X)$.

\proclaim{Theorem 4.6} Let  $\phi\in  L^{\infty} (\jj, X)$ and $u$ be a bounded mild solution of (4.12). Let $\A $ satisfy (3.1), (3.6), $\gamma_{\la}\A \st \A$ for all $\la\in \r$ and contain all constants.

(i) If $\jj=\r_+$, then $i\, sp^{\Cal {L}} (u)\st (\sigma(A)\cap i\r)\cup i\,sp^{\Cal {L}}(\phi)$.

(ii) If $sp_{\A} (\phi)=\emptyset$, then
  $i\,sp_{\A} (u) \st \sigma(A)\cap i\r$.
\endproclaim

\demo{Proof}  As $u, \phi\in L^{\infty} (\jj,X)$ we get  $M_h u, M_h \phi\in BUC(\jj,X)$  and $v= M_h u$ is a classical solution of $ v'(t)= A v(t) + M_h  \phi (t) $, $v(t) \in D(A)$, $t\in  {\jj}$ for each $h >0$.

(i)  By [5, Proposition 5.6.7, p. 380], we have

$i\,sp^{\Cal {L}} (M_h u)\st (\sigma(A)\cap i\r)\cup i\,sp^{\Cal {L}} (M_h\phi)$ for all $h >0$.

\noindent Taking the union of both sides, we get

 $\cup _{h >0}\, i\,sp^{\Cal {L}} (M_h u)\st (\sigma(A)\cap i\r)\cup (\cup _{h >0} \, i\,sp^{\Cal {L}} (M_h (\phi))$.

\noindent Applying Proposition 4.1 (ii) to both sides, we conclude that

$i\,sp^{\Cal {L}} (u)\st (\sigma(A)\cap i\r)\cup i\,sp^{\Cal {L}} (\phi)$.

 (ii)  Take $h >0$. Since $sp_{\A} (\phi)=\emptyset$, it follows that $sp_{\A} (M_h\phi)=\emptyset$,  by Proposition 3.4 (ii). Hence $M_h\phi\in \A$ by [6, Theorem 4.2.1].  Using [7, Corollary 3.4 (i)], we conclude that  $i\,sp_{\A} (M_h u) \st \sigma(A)\cap i\r$. By Proposition 3.4(ii), we conclude that $i\,sp_{\A} (u) \st \sigma(A)\cap i\r$.
  \P
\enddemo

\proclaim{Remark 4.7} (i) There are  some inclusions of this general type in [4, (4.4), (4.5)].

(ii) In [11], the inclusion
  $i\,sp_{\A} (u) \st (\sigma(A)\cap i\r) \cup i\,sp_{\A} (\phi)$ was proved.
\endproclaim

\smallskip

\smallskip

\indent School of Math. Sci., P.O. Box 28M, Monash University,
 Vic. 3800.

\indent Email "bolis.basit\@monash.edu",\qquad
"alan.pryde\@monash.edu".

\enddocument